\documentclass[11pt,a4paper]{article}
\usepackage[T1]{fontenc}
\usepackage[utf8]{inputenc}
\usepackage[tbtags]{amsmath}
\usepackage{amsthm,amsfonts,amssymb,mathrsfs,amscd}
\usepackage[margin=28mm]{geometry}
\usepackage[hidelinks]{hyperref}
\numberwithin{equation}{section}
\theoremstyle{plain}
\newtheorem{theorem}{Theorem}
\newtheorem{lemma}{Lemma}[section]

\newtheorem{corollary}{Corollary}
\newtheorem{thmx}{Theorem}

\theoremstyle{definition}

\theoremstyle{remark}
\newtheorem{remark}{Remark}

\begin{document}
\title{On the relationship between two classical problems concerning the almost-everywhere behavior of Fourier sums}
\author{D.\,I.~Masyutin\\[3pt]
\small N. N. Krasovskii Institute of Mathematics and Mechanics,\\
\small Ural Branch of the Russian Academy of Sciences, Yekaterinburg, Russia\\
\small \texttt{newselin@mail.ru}}
\date{November 26, 2026}
\maketitle
\begin{abstract}
The paper deals with two classical problems in the theory of trigonometric Fourier series: Hardy's problem concerning the exact order of growth of partial Fourier sums of integrable functions on a set of full measure, and Zygmund's problem concerning sharp sufficient conditions on the integral modulus of continuity that ensure the almost everywhere convergence of Fourier series. We establish a relationship between Zygmund's problem and a certain condition on partial Fourier sums, which in turn is closely related to Hardy's problem. The result obtained refines the well-known theorems of Tandori and Totik.
\end{abstract}


\section{Introduction}
\subsection{Notation and definitions}
Let $\mathbb{T} = [-\pi, \pi)$ be the one-dimensional torus, and let $\mathbb{N}$ and $\mathbb{Z}_+$ denote the sets of positive and nonnegative integers, respectively. Throughout the paper, $\psi \colon [0, +\infty) \to (0, +\infty)$ denotes a nondecreasing function, and $\chi \colon [0, +\infty) \to (0, +\infty)$ denotes an arbitrary positive function. We define the class $L\psi(L)(\mathbb{T})$ as the set of Lebesgue measurable real-valued functions $f$ on the torus satisfying
\begin{equation*}
\int_{\mathbb{T}} |f(x)| \psi(|f(x)|) \, dx < \infty.
\end{equation*}

For a function $f \in L(\mathbb{T})$, consider its trigonometric Fourier series
\begin{equation}\label{fourier}
\frac{a_0}{2} + \sum_{k=1}^{\infty} (a_k \cos{kx} + b_k \sin{kx}),
\end{equation}
where the coefficients $a_k$ and $b_k$ are given by the standard formulas
\begin{equation*}
a_0 = \frac{1}{\pi} \int\limits_{\mathbb{T}} f(x) \, dx,\quad a_k = \frac{1}{\pi} \int\limits_{\mathbb{T}} f(x) \cos{kx}\, dx,\quad b_k = \frac{1}{\pi} \int\limits_{\mathbb{T}} f(x) \sin{kx}\, dx,\quad k \in \mathbb{N}.
\end{equation*}
We denote by $S_n(f, x)$ the value of the $n$th partial sum of the series~\eqref{fourier} at $x \in \mathbb{T}$:
\begin{equation*}
S_n(f, x) = \frac{a_0}{2} + \sum_{k=1}^{n} (a_k \cos{kx} + b_k \sin{kx}),\qquad n \in \mathbb{Z}_+.
\end{equation*}

We write $mE$ for the Lebesgue measure of a measurable set $E \subset \mathbb{T}$. We also introduce the following notation:
\begin{gather*}
\|f\|_1 := \int_{\mathbb{T}} |f(x)| \, dx,\qquad f^+(x) := \max\{f(x), 0\}, \\
M_n(\psi, \chi)(f, x) := \frac{1}{\psi(n)} \max_{0 \leqslant m \leqslant n} \frac{|S_m(f, x)|}{\chi(m)},\qquad
T^*(f, x) := \sup_{n \in \mathbb{Z}_+} |T_n(f, x)|,
\end{gather*}
where $\{T_n \colon L(\mathbb{T}) \to L(\mathbb{T})\}_{n \in \mathbb{Z}_+}$ is an arbitrary sequence of operators. Recall that the weak-type $(1,1)$ inequality for an operator
\(T\colon L(\mathbb T)\to L(\mathbb T)\) takes the form
\[
m\{x\in\mathbb T:\ |T(f,x)|>y\}\le \frac{C}{y}\|f\|_1
\]
for all \(f\in L(\mathbb T)\) and \(y>0\), where \(C>0\) is a constant
independent of \(f\) and \(y\). We say that a sequence of operators $\{T_n\}_{n \in \mathbb{Z}_+}$ satisfies a weak-type $(1,1)$ inequality if there exists a constant $C>0$ such that, for all $f \in L(\mathbb{T})$, $n \in \mathbb{Z}_+$, and $y>0$,
\begin{equation*}
m\bigl\{x \in \mathbb{T}: |T_n(f, x)| > y\bigr\} \leqslant \frac{C}{y} \|f\|_1.
\end{equation*}

Next, for a function $f \in L(\mathbb{T})$, define its modulus of continuity in $L(\mathbb{T})$ by
\begin{equation*}
\omega(f, \delta)_1 = \sup_{h\in\mathbb{R}:\, |h| \leqslant \delta} \| f(\cdot+h) - f \|_1,
\end{equation*}
and its best approximation by trigonometric polynomials of degree at most $n$ by
\begin{equation*}
E_n(f)_1 = \inf_{T_n \in \mathfrak{T}_n} \|f - T_n\|_1,
\end{equation*}
where $\mathfrak{T}_n$ is the space of all trigonometric polynomials of degree at most $n$. By a modulus of continuity $\omega$ we mean a subadditive, nondecreasing, nonnegative function on $[0, +\infty)$ satisfying
\begin{equation*}
\omega(0) = 0 = \lim\limits_{\delta \to +0} \omega(\delta).
\end{equation*}
For such a function $\omega$, define the class
\begin{equation*}
H_1^{\omega} := \bigl\{f \in L(\mathbb{T}): \omega(f, \delta)_1 = O(\omega(\delta))\bigr\}.
\end{equation*}

\subsection{Some classical results on the almost everywhere behavior of Fourier sums}

In 1913, G.\,H.~Hardy~\cite{Hardy} obtained the following estimate for the growth of Fourier sums: for every $f \in L(\mathbb{T})$ and almost every $x \in \mathbb{T}$,
\begin{equation}\label{est_log}
S_n(f, x) = o(\ln{n}),\quad n \to \infty.
\end{equation}
In the same paper, Hardy conjectured that this estimate is sharp. In 1985, K.\,I.~Oskolkov~\cite{Oskolkov} generalized~\eqref{est_log} by proving that, for every increasing sequence of positive integers $\{n_k\}_{k=1}^{\infty}$, the relation $S_{n_k}(f, x) = o(\ln{k})$ as $k \to \infty$ holds almost everywhere. On the other hand, A.\,N.~Kolmogorov proved that there exist functions in $L(\mathbb{T})$ whose Fourier series diverge almost everywhere~\cite{Kolmogoroff1}, and even everywhere~\cite{Kolmogoroff2}. In 1962, Y.\,M.~Chen~\cite{Chen1} made progress on Hardy's conjecture by proving that, for every sequence $\psi(n) = o(\ln{\ln{n}})$ as $n \to \infty$, there exists a function $F \in L(\mathbb{T})$ such that, for every $x \in \mathbb{T}$, the inequality $|S_n(F, x)| > \psi(n)$ holds for infinitely many values of $n$. In the same paper, Chen conjectured the existence of a function in $L\ln^+{\ln^+{(L)}}(\mathbb{T})$ with an everywhere divergent Fourier series. Subsequently, in 1969, he~\cite{Chen2} (see also~\cite{Prohorenko}) constructed a function in $L(\ln^+{\ln^+{(L)}})^{1-\varepsilon}(\mathbb{T})$, $0 < \varepsilon \leqslant 1$, whose Fourier series diverges almost everywhere. In 1981, T.\,W.~K\"orner~\cite{Korner} strengthened this result somewhat: for every function $\psi(n) = o(\ln{\ln{n}})$ as $n \to \infty$, there exists a function in $L\psi(L)(\mathbb{T})$ with an everywhere divergent Fourier series. The strongest result in this direction to date is due to S.\,V.~Konyagin~\cite{Konyagin}, who proved that, for any nondecreasing functions $\psi$ and $\chi$ such that $\psi(n)\chi(n) = o(\sqrt{\ln{n}/\ln{\ln{n}}})$ as $n \to \infty$, there exists a function $F \in L\psi(L)(\mathbb{T})$ satisfying
\begin{equation}\label{not_est}
\varlimsup\limits_{n \to \infty} \frac{|S_n(F, x)|}{\chi(n)} = +\infty
\end{equation}
for every $x \in \mathbb{T}$. The question of the sharpness of estimate~\eqref{est_log} remains open.

In his monograph, A.~Zygmund~\cite[Ch.~13, Theorem~3.10]{Zigmund} proved that if, for a modulus of continuity $\omega$, the integral
$
\int\limits_0^1 \omega(t)t^{-1} \, dt
$
converges, or equivalently,
\begin{equation}\label{cond_1}
\sum\limits_{n=1}^{\infty} \frac{\omega(1/n)}{n} < +\infty,
\end{equation}
then the Fourier series of every $f \in H_1^{\omega}$ converges almost everywhere. He also raised the question of whether condition~\eqref{cond_1} can be weakened. In the aforementioned paper~\cite{Oskolkov}, K.\,I.~Oskolkov proved, in particular, that condition~\eqref{cond_1} can be replaced by
\begin{equation}\label{cond_2}
\sum\limits_{n=1}^{\infty} \frac{E_n(f)_1}{n} < +\infty.
\end{equation}
Moreover, the derivation of this condition from an estimate of the form~\eqref{est_log} can be traced in~\cite{Oskolkov}.

In connection with the sharpness of condition~\eqref{cond_1}, V.\,I.~Prokhorenko~\cite{Prohorenko} constructed in 1968 a function in the class $H_1^{\overline{\omega}}$, where $\overline{\omega}(\delta) = 1/\ln{\ln{(1/\delta)}}$, whose Fourier series diverges almost everywhere. As Prokhorenko observed, by an embedding theorem of P.\,L.~Ul'yanov~\cite[Theorem~1]{Ulyanov}, this result implies the aforementioned 1969 result of Chen. In 2004, N.\,Yu.~Antonov~\cite{Antonov}, using the main lemma of~\cite{Konyagin}, constructed a function in the class $H_1^{\omega^*}$, where $\omega^*(\delta) = \sqrt{\ln{\ln{(1/\delta)}}/\ln{(1/\delta)}}$, with an almost everywhere divergent Fourier series. Nevertheless, the question of the sharpness of conditions~\eqref{cond_1} and~\eqref{cond_2} also remains open.

\subsection{The relationship between the problems of the almost everywhere behavior of Fourier sums}

The results above naturally suggest that progress on either of the two problems under consideration should lead to progress on the other. In 1983, K.~Tandori~\cite{Tandori} proved that if, for a sufficiently regular function $\psi$, there exists a function $F \in L(\mathbb{T})$ such that
\begin{equation}\label{Tandori_eq}
\varlimsup\limits_{n \to \infty} \frac{|S_n(F, x)|}{\psi(n)} > 0
\end{equation}
on a set of positive measure, then there exists a function in $L\psi^{1-\varepsilon}(L)(\mathbb{T})$, where $0 < \varepsilon \leqslant 1$, whose Fourier series diverges almost everywhere. In the same year, V.~Totik~\cite{Totik} complemented this result by showing that the hypotheses of Tandori's theorem imply the existence of a function in $H_1^{1/\psi}$ with an almost everywhere divergent Fourier series. Moreover, he established a partial converse: if, for some $\varepsilon > 0$, there exists a function in $H_1^{1/\psi^{1+\varepsilon}}$ whose Fourier series diverges on a set of positive measure, then there exists a function in $L(\mathbb{T})$ satisfying~\eqref{Tandori_eq} almost everywhere.

The main purpose of the present paper is to prove the following theorem, which refines the results of Tandori and Totik.

\begin{theorem}\label{theorem1}
Let $\psi, \chi \colon [0, +\infty) \to (0, +\infty)$ be nondecreasing functions such that the sequence $n(\psi(n+1) - \psi(n))$ is nonincreasing. Then the following conditions are equivalent:
\begin{enumerate}
\item\label{item:first} The sequence of operators $\{M_n(\psi, \chi)\}_{n \in \mathbb{Z}_+}$ satisfies a weak-type $(1, 1)$ inequality.
\item\label{item:second} For every function $f \in L(\mathbb{T})$ satisfying
\begin{equation}\label{series_cond}
\sum\limits_{n=1}^{\infty} [\psi(n+1) - \psi(n)] \omega\left(f, \frac{1}{n}\right)_1 < +\infty,
\end{equation}
the following relation holds for almost every $x \in \mathbb{T}$:
\begin{equation}\label{max_S_n}
\lim\limits_{n \to \infty} \frac{|f(x) - S_n(f, x)|}{\chi(n)} = 0.
\end{equation}
\end{enumerate}
\end{theorem}
\begin{remark}
Let $\psi(n)\chi(n)=\ln n$. If $\chi\equiv c>0$, then $\psi(n)=c^{-1}\ln n$, and~\eqref{series_cond} coincides with Zygmund's condition $\sum_{n=1}^{\infty} \omega(f,1/n)_1/n <\infty$, while~\eqref{max_S_n} is equivalent to the almost everywhere convergence of the Fourier series. If, on the other hand, $\psi\equiv c>0$, then $\chi(n)=c^{-1}\ln n$, and condition~\ref{item:first} is a weak-type $(1,1)$ inequality for $\sup |S_m(f,x)|/\ln m$, which implies Hardy's estimate $S_n(f,x)=o(\ln n)$ almost everywhere by a standard argument (see, e.g.,~\cite[Ch.~2, \S~3, item~3.2]{Torchinsky}). Thus, both classical problems are encompassed by Theorem~\ref{theorem1}.
\end{remark}
\begin{remark}
By the aforementioned theorem of Ul'yanov~\cite[Theorem~1]{Ulyanov}, condition~\eqref{series_cond} implies that $f \in L\psi(L)(\mathbb{T})$. It follows that if the sequence of operators $\{M_n(\psi, \chi)\}_{n \in \mathbb{Z}_+}$ fails to satisfy a weak-type $(1, 1)$ inequality, then there exists a function $f \in L\psi(L)(\mathbb{T})$ for which~\eqref{not_est} holds on a set of positive measure.
\end{remark}

In the final part of the paper, we derive the following result from the results obtained here and the key lemma of~\cite{Konyagin}.
\begin{corollary}\label{theorem2}
Let $\psi, \chi: [0, +\infty) \rightarrow (0, +\infty)$ be nondecreasing functions such that $\psi(n)\chi(n) = o(\sqrt{\ln{n}/\ln{\ln{n}}})$ as $n \to \infty$ and $n(\psi(n+1) - \psi(n))$ is nonincreasing on $\mathbb{N}$. Then there exists a function $F \in L(\mathbb{T})$ satisfying~\eqref{series_cond} for which~\eqref{not_est} holds almost everywhere.
\end{corollary}
\medskip

\section{Preliminaries}

Before proceeding to the proofs, we explain how Theorem~\ref{theorem1} is related to the well-known theorem of E.~Stein~\cite{Stein}. We state this result in the particular form needed here. Recall that an operator $T \colon L(\mathbb{T}) \to L(\mathbb{T})$ is called translation invariant if
\begin{equation*}
T(f(\cdot + s), x) := T(f(\cdot), x+s)
\end{equation*}
for all $s \in \mathbb{T}$.
\begin{thmx}[{\cite[Theorem~1]{Stein}}]\label{t_Stein}
Let $\{T_n\}_{n \in \mathbb{Z}_+}$ be a sequence of bounded linear translation-invariant operators from $L(\mathbb{T})$ to $L(\mathbb{T})$. If, for every function $f \in L(\mathbb{T})$ and almost every $x \in \mathbb{T}$, the sequence $T_n(f, x)$ converges, then the operator $T^*(f, x)$ satisfies a weak-type $(1, 1)$ inequality: there exists a constant $C>0$ such that, for every $f \in L(\mathbb{T})$ and every $y>0$,
\begin{equation*}
m\left\{x \in \mathbb{T}: T^*(f, x) > y \right\} \leqslant \frac{C}{y} \|f\|_1.
\end{equation*}
\end{thmx}
Applying Theorem~\ref{t_Stein} to the sequence of linear operators
\begin{equation*}
T_n(f, x) := \frac{f(x) - S_n(f, x)}{\chi(n)},
\end{equation*}
we immediately obtain one of the implications in Theorem~\ref{theorem1} in the case where $\psi$ is constant. Thus, Theorem~\ref{theorem1} is an analogue of Stein's theorem for functions in the spaces $H_1^{\omega}$.

Note that a standard argument (see, e.g.,~\cite[Ch.~2, \S~3, item~3.2]{Torchinsky}) shows that if $\chi(n) \to \infty$ as $n\to\infty$ and the expression
\begin{equation}\label{majorant}
\sup\limits_{n \in \mathbb{Z}_+} \frac{|S_n(f, x)|}{\chi(n)},
\end{equation}
is finite almost everywhere, then for every $f \in L(\mathbb{T})$ the following relation holds almost everywhere:
\begin{equation*}
S_n(f, x) = o(\chi(n)),\quad n \to \infty.
\end{equation*}
In particular, if $\chi(n) = \ln{n}$, the finiteness of~\eqref{majorant} on a set of full measure is equivalent to Hardy's relation~\eqref{est_log} on a set of full measure. Note also that if~\eqref{majorant} is finite almost everywhere, then, by the preceding observation and Theorem~\ref{t_Stein}, the operator~\eqref{majorant} satisfies a weak-type $(1, 1)$ inequality. Hence, for any functions $\psi$ and $\theta$ such that $\psi(n)\theta(n) = \chi(n)$, the sequence of operators $\{M_n(\psi, \theta)\}_{n \in \mathbb{Z}_+}$ also satisfies a weak-type $(1, 1)$ inequality. In particular, it follows that a sufficient condition of the form~\eqref{cond_1} for the almost everywhere convergence of Fourier series follows from a relation of the form~\eqref{est_log}. The converse question, however, remains open: whether a sufficient condition of the form~\eqref{cond_1} for the almost everywhere convergence of Fourier series implies a relation of the form~\eqref{est_log}, or, equivalently, whether the weak-type $(1,1)$ inequality for the sequence of operators $\{M_n(\psi, \theta)\}_{n \in \mathbb{Z}_+}$ implies the almost everywhere finiteness of the operator~\eqref{majorant}.

\subsection{A lemma on series}
The following result shows, in particular, that conditions~\eqref{cond_1} and~\eqref{cond_2} are equivalent.

\begin{lemma}\label{lemma_series}
Let $f \in L(\mathbb{T})$, and let $\alpha: [0, +\infty) \rightarrow [0, +\infty)$ be a function such that $n^{p}\alpha(n)$ is nonincreasing for some $p>0$. Then the series
\begin{align}
\sum\limits_{n=1}^{\infty} \alpha(n) \omega\left(f, \frac{1}{n}\right)_1,\label{series1} \\
\sum\limits_{n=1}^{\infty} \alpha(n) E_n(f)_1 \label{series2}
\end{align}
either both converge or both diverge.
\end{lemma}

\par\noindent\textit{Proof.}\enspace
Recall the well-known Jackson inequality in $L(\mathbb{T})$~\cite{Quade} (see also~\cite[Appendices, \S~7]{Bari}): for every function $f \in L(\mathbb{T})$ there exists a constant $C>0$ such that
\begin{equation*}
E_n(f)_1 \leqslant C \omega\left(f, \frac{1}{n}\right)_1,\quad n \in \mathbb{N}.
\end{equation*}
Consequently, the convergence of~\eqref{series1} implies the convergence of~\eqref{series2}, and the divergence of~\eqref{series2} implies the divergence of~\eqref{series1}.

Next, we use the inequality of A.\,F.~Timan and M.\,F.~Timan~\cite[Ch.~6, item~6.1.1]{Timan}: there exists a constant $A > 0$ such that
\begin{equation}\label{Timan_eq}
\omega\left(f, \frac{1}{n}\right)_1 \leqslant \frac{A}{n} \sum\limits_{k=1}^n E_k(f)_1,\quad n \in \mathbb{N}.
\end{equation}

Suppose that the series~\eqref{series2} converges. Using inequality~\eqref{Timan_eq} and the fact that $n^p \alpha(n)$ is nonincreasing, we obtain
\begin{equation}\label{omega_E_n}
\sum\limits_{n=1}^{\infty} \alpha(n) \omega\left(f, \frac{1}{n}\right)_1 \leqslant A \sum\limits_{k=1}^{\infty} E_k(f)_1 \sum\limits_{n=k}^{\infty} \frac{\alpha(n)}{n}.
\end{equation}
Since $n^p \alpha(n)$ is nonincreasing, we have
\begin{equation*}
\sum\limits_{n=k}^{\infty} \frac{\alpha(n)}{n} \leqslant k^p \alpha(k) \sum\limits_{n=k}^{\infty} \frac{1}{n^{p+1}} \leqslant \alpha(k) \left(\frac{1}{p} + \frac{1}{k}\right).
\end{equation*}
This and estimate~\eqref{omega_E_n} yield
\begin{equation*}
\sum\limits_{n=1}^{\infty} \alpha(n) \omega\left(f, \frac{1}{n}\right)_1 \leqslant A \frac{p+1}{p} \sum\limits_{k=1}^{\infty} \alpha(k) E_k(f)_1.
\end{equation*}
Hence, the series~\eqref{series1} converges. This proves Lemma~\ref{lemma_series}.

\begin{remark}\label{remark3}
Let $\psi \colon [0, +\infty) \to (0, +\infty)$ be a nondecreasing function such that the sequence $n(\psi(n+1) - \psi(n))$ is nonincreasing. Then, by Lemma~\ref{lemma_series}, condition~\eqref{series_cond} is equivalent to the following condition:
\begin{equation}\label{series_E_n(f)}
\sum\limits_{n=1}^{\infty} [\psi(n+1) - \psi(n)] E_n(f)_1 < +\infty.
\end{equation}
\end{remark}

\subsection{Auxiliary results}

The proof of the implication from condition~\ref{item:second} to condition~\ref{item:first} in Theorem~\ref{theorem1} is based on the proof of Theorem~\ref{t_Stein}. It consists of two steps, which we formulate as lemmas. In proving the first lemma, we use the following well-known probabilistic result~\cite[Lemma~19]{Paley}.

\begin{thmx}\label{anticoncentration}
Let $t \in [0, 1]$, $n \in \mathbb{N}$, let $r_n(t):= \operatorname{sign}({\sin{(2^n\pi t)}})$ be the Rademacher functions, and let $\mu$ denote Lebesgue measure on $[0, 1]$. Then, for every sequence of real numbers $\{x_n\}_{n=1}^{\infty}$ satisfying
\begin{equation*}
\sum\limits_{n=1}^{\infty} x_n^2 < \infty,
\end{equation*}
the following inequality holds:
\begin{equation*}
\mu\left\{t \in [0, 1]: \left|\sum\limits_{n=1}^{\infty} x_n r_n(t) \right| > \frac{1}{2} \left(\sum\limits_{n=1}^{\infty} x_n^2 \right)^{1/2} \right\} \geqslant \frac{1}{16}.
\end{equation*}
\end{thmx}  

\begin{lemma}\label{lemma}
Let $\psi$ be a positive nondecreasing function on $[0, +\infty)$, and let $\chi$ be an arbitrary positive function on $[0, +\infty)$. Suppose that the sequence of operators $\{M_n(\psi, \chi)\}_{n=1}^{\infty}$ fails to satisfy a weak-type $(1, 1)$ inequality. Then there exist a sequence of positive integers $\{p_n\}_{n=1}^{\infty}$, a sequence of real numbers $\{A_n\}_{n=1}^{\infty}$ tending to infinity, a sequence of sets $\{E_n\}_{n=1}^{\infty}$, and a sequence of trigonometric polynomials $\{f_n\}_{n=1}^{\infty}$ satisfying the following conditions:
\begin{enumerate}
\item\label{item0} There exists $c>0$ such that $mE_n \geqslant c$ for every $n \in \mathbb{N}$.
\item\label{item1} The degree of each $f_n$ is at most $2p_n$.
\item\label{item2} There exists $C>0$ such that $\|f_n\|_1 \leqslant C$ for every $n \in \mathbb{N}$.
\item For every $x \in E_n$,
\begin{equation}\label{S>m}
\max\limits_{0 \leqslant m \leqslant p_n} \frac{|S_m(f_n, x)|}{\chi(m)} > \psi(p_n) A_n.
\end{equation}
\end{enumerate}
\end{lemma}

\par\noindent\textit{Proof.}\enspace
The hypothesis implies that for every $C > 0$ there exist a function $F_C \in L(\mathbb{T})$ with $\|F_C\|_1 = 1$, a positive integer $n = n(C) \in \mathbb{N}$, and a number $y(C) > 0$ such that
\begin{equation*}
m\left\{x \in \mathbb{T}: \max\limits_{0 \leqslant m \leqslant n} \frac{|S_m(F_C, x)|}{\chi(m)} > \psi(n) y(C) \right\} \geqslant \frac{C}{y(C)}. 
\end{equation*}
Setting $C = 2^n$, $n \in \mathbb{N}$, we obtain a sequence of positive integers $\{p_n\}_{n=1}^{\infty}$ such that the set
\begin{equation}\label{G_n_def}
G_n := \left\{x\in \mathbb{T}: \max\limits_{0 \leqslant m \leqslant p_n} \frac{|S_m(F_{2^n}, x)|}{\chi(m)} > \psi(p_n) y(2^n) \right\}
\end{equation}
satisfies the estimate
\begin{equation}\label{m_not_weak}
\frac{2^n}{y(2^n)} \leqslant mG_n. 
\end{equation}

Define a sequence of polynomials $\{g_n\}_{n=1}^{\infty}$ by setting
\begin{equation}\label{f_k_def}
g_n(x) := \sigma_{2p_n, p_n}(F_{2^n}, x) := \frac{1}{p_n+1} \sum_{m=p_n}^{2p_n} S_m(F_{2^n}, x),
\end{equation}
where $\sigma_{2p_n, p_n}$ is the de la Vall\'ee Poussin mean. It can be represented in terms of Fej\'er means as follows:
\begin{equation}\label{representation}
\sigma_{2p_n, p_n}(F_{2^n}, x) = \frac{2p_n + 1}{p_n+1} \sigma_{2p_n}(F_{2^n}, x) - \frac{p_n}{p_n+1} \sigma_{p_n-1}(F_{2^n}, x).
\end{equation}
From~\eqref{representation} and the estimate for the norm of the Fej\'er kernel, it follows that
\begin{equation}\label{norm_g}
\|g_n\|_1 \leqslant 3.
\end{equation}
Note also that, by definition~\eqref{f_k_def}, for every $x \in \mathbb{T}$ and every $m \leqslant p_n$ we have
\begin{equation}\label{eq_S}
S_m(g_n, x) = S_m(F_{2^n}, x).
\end{equation}

Construct a sequence $\{s_n\}_{n=1}^{\infty}$ as follows:
\begin{equation}\label{s_k}
s_n := \min\left\{m \in \mathbb{N}: 1 \leqslant \frac{2^n}{y(2^n)} m \leqslant 2\pi \right\}.
\end{equation}
It follows immediately from~\eqref{s_k} that
\begin{gather}
\frac{2^n}{y(2^n)} s_n \geqslant 1, \label{sum_E} \\
\frac{y(2^n)}{s_n} \geqslant \frac{2^{n-1}}{\pi}. \label{sum_M_k}
\end{gather}

Next, for a fixed tuple $(x_1, \ldots, x_{s_n}) \in \mathbb{T}^{s_n}$, define
\begin{equation*}
F(x_1, \ldots, x_{s_n}) := m\left(\bigcup\limits_{j=1}^{s_n} G_n^{x_j}\right) = \int\limits_{\mathbb{T}} I_{\bigcup\limits_{j=1}^{s_n} G_n^{x_j}}(y)\, dy, 
\end{equation*}
where $I_M$ denotes the characteristic function of the set $M$ and $G_n^{x_j} := G_n + x_j$ is the translate of $G_n$ by $x_j$. Representing $F$ as an integral of a characteristic function, then integrating with respect to $(x_1, \ldots, x_{s_n})$ and interchanging the order of integration, we obtain
\begin{equation}\label{int_F}
\frac{1}{(2\pi)^{s_n}}\int\limits_{\mathbb{T}^{s_n}} F(x_1, \ldots, x_{s_n}) \, dx_1 \ldots dx_{s_n} = \int\limits_{\mathbb{T}} \frac{1}{(2\pi)^{s_n}}\int\limits_{\mathbb{T}^{s_n}} I_{\bigcup\limits_{j=1}^{s_n} G_n^{x_j}}(y) \, dx_1 \ldots dx_{s_n} dy.
\end{equation}
Using the relation
\begin{equation*}
I_{\bigcup\limits_{j=1}^{s_n} G_n^{x_j}}(y) = 1 - \prod\limits_{j=1}^{s_n} I_{\mathbb{T} \setminus G_n^{x_j}}(y)
\end{equation*}
and identity~\eqref{int_F}, we obtain
\begin{equation}\label{int_F2}
\frac{1}{(2\pi)^{s_n}}\int\limits_{\mathbb{T}^{s_n}} F(x_1, \ldots, x_{s_n}) \, dx_1 \ldots dx_{s_n} = 2\pi - 2\pi\left(1 - \frac{mG_n}{2\pi}\right)^{s_n}.
\end{equation}
Identity~\eqref{int_F2} and inequalities~\eqref{m_not_weak} and~\eqref{sum_E} imply that
\begin{equation*}
\frac{1}{(2\pi)^{s_n}}\int\limits_{\mathbb{T}^{s_n}} F(x_1, \ldots, x_{s_n}) \, dx_1 \ldots dx_{s_n} \geqslant 2\pi - 2\pi\left(1 - \frac{2^n}{2\pi y(2^n)}\right)^{y(2^n)/2^n} \geqslant 2\pi (1 - e^{-1/(2\pi)}).
\end{equation*}
Consequently, there exist a constant $c_1>0$ and a tuple $(\overline{x}_1, \ldots, \overline{x}_{s_n}) \in \mathbb{T}^{s_n}$ such that
\begin{equation}\label{ineq_meas}
m\left(\bigcup\limits_{j=1}^{s_n} G_n^{\overline{x}_j}\right) \geqslant c_1.
\end{equation}

Let
\begin{equation}\label{f_n^t_def}
f^n_t(x) := \frac{1}{s_n} \sum\limits_{j=1}^{s_n} r_j(t) g_n(x-\overline{x}_j),\quad t \in [0, 1],
\end{equation}
where $r_j(t):= \operatorname{sign}({\sin{(2^j\pi t)}})$ are the Rademacher functions. Suppose that $x \in G_n^{\overline{x}_j}$ for some $j \in \{1, \ldots, s_n\}$. Since $x - \overline{x}_j \in G_n$, it follows from~\eqref{eq_S} and~\eqref{G_n_def} that there exists $m(x) \leqslant p_n$ such that
\begin{equation}\label{S>m2}
\frac{|S_{m(x)}(g_n, x - \overline{x}_j)|}{\chi(m(x))} > \psi(p_n) y(2^n).
\end{equation}
Apply Theorem~\ref{anticoncentration} to a Fourier partial sum of the function~\eqref{f_n^t_def}. Then
\begin{equation*}
\mu\left\{t \in [0, 1]: \frac{\left|S_{m(x)}(f^n_t, x)\right|}{\chi(m(x))} > \frac{1}{2 s_n} \left(\sum\limits_{j=1}^{s_n} \left(\frac{S_{m(x)}(g_n, x-\overline{x}_j)}{\chi(m(x))}\right)^2 \right)^{1/2} \right\} \geqslant \frac{1}{16}.
\end{equation*}
Thus, by inequality~\eqref{S>m2},
\begin{equation}\label{anticon_ineq}
\mu\left\{t \in [0, 1]: \max\limits_{0 \leqslant m \leqslant p_n} \frac{\left|S_m(f^n_t, x) \right|}{\chi(m)} > \frac{y(2^n)}{2 s_n} \psi(p_n) \right\} \geqslant \frac{1}{16},\quad x \in \bigcup\limits_{j=1}^{s_n} G_n^{\overline{x}_j}.
\end{equation}
Again, writing the right-hand side of inequality~\eqref{anticon_ineq} as an integral of a characteristic function, integrating both sides of~\eqref{anticon_ineq} over $\bigcup_{j=1}^{s_n} G_n^{\overline{x}_j}$, interchanging the order of integration, and applying inequality~\eqref{ineq_meas}, we obtain
\begin{multline}\label{anticon_ineq_2}
\int\limits_{[0, 1]} m\left\{x \in \bigcup_{j=1}^{s_n} G_n^{\overline{x}_j}: \max\limits_{0 \leqslant m \leqslant p_n} \frac{\left|S_m(f^n_t, x) \right|}{\chi(m)} > \frac{y(2^n)}{2 s_n} \psi(p_n) \right\}\, dt \\
\geqslant \frac{1}{16} m\left(\bigcup\limits_{j=1}^{s_n} G_n^{\overline{x}_j} \right) \geqslant \frac{c_1}{16}.
\end{multline}
Hence, there exist a constant $c>0$ and a point $\overline{t} \in [0, 1]$ such that
\begin{equation}\label{anticon_ineq_3}
m\left\{x \in \bigcup\limits_{j=1}^{s_n} G_n^{\overline{x}_j}: \max\limits_{0 \leqslant m \leqslant p_n} \frac{\left|S_m(f^n_{\overline{t}}, x) \right|}{\chi(m)} > \frac{y(2^n)}{2 s_n} \psi(p_n) \right\} \geqslant c.
\end{equation}

Set
\begin{gather*}
f_n(x) := f_{\overline{t}}^n(x), \\
E_n :=  \left\{x \in \bigcup\limits_{j=1}^{s_n} G_n^{\overline{x}_j}: \max\limits_{0 \leqslant m \leqslant p_n} \frac{\left|S_m(f^n_{\overline{t}}, x) \right|}{\chi(m)} > \frac{y(2^n)}{2 s_n} \psi(p_n) \right\}.
\end{gather*}
By~\eqref{anticon_ineq_3}, condition~\ref{item0} holds. Inequalities~\eqref{anticon_ineq_2} and~\eqref{sum_M_k} imply that, for $x \in E_n$,
\begin{equation}\label{ineq_3}
\max\limits_{0 \leqslant m \leqslant p_n} \frac{\left|S_m(f^n_{\overline{t}}, x) \right|}{\chi(m)} > \frac{2^n}{4 \pi} \psi(p_n).
\end{equation}
Thus, inequality~\eqref{S>m} follows from~\eqref{ineq_3}. Recalling definitions~\eqref{f_n^t_def} and~\eqref{f_k_def}, we immediately obtain condition~\ref{item1}. Condition~\ref{item2}, in turn, follows from definition~\eqref{f_n^t_def} and inequality~\eqref{norm_g}. This proves Lemma~\ref{lemma}.

To prove the second lemma, we need the following classical result~\cite[Ch.~13, \S~1, Lemma~1.24]{Zigmund}.
\begin{thmx}\label{Borel-Kantelli}
Let $\{ E_n \}_{n=1}^{\infty}$ be a sequence of measurable subsets of $\mathbb{T}$ such that
\begin{equation}\label{m_sum}
\sum\limits_{n=1}^{\infty} mE_n = \infty.
\end{equation}
Then there exists a sequence $\{ x_n \}_{n=1}^{\infty} \subset \mathbb{T}$ with the property
\begin{equation*}
m\left(\varlimsup\limits_{n \to \infty} E_n^{x_n} \right) := m\left(\bigcap\limits_{k=1}^{\infty} \bigcup\limits_{n=k}^{\infty} E_n^{x_n} \right) = 2 \pi,
\end{equation*}
where $E_n^{x_n} := \{x+x_n: x\in E_n \}$ denotes the translate of $E_n$ by $x_n$.
\end{thmx}

\begin{lemma}\label{lemma2}
Suppose that the conclusion of Lemma~\ref{lemma} holds with positive nondecreasing functions $\psi$ and $\chi$ such that the sequence $n (\psi(n+1) - \psi(n))$ is nonincreasing. Then there exists a function $F \in L(\mathbb{T})$ satisfying condition~\eqref{series_cond} such that relation~\eqref{not_est} holds almost everywhere.
\end{lemma}

\par\noindent\textit{Proof.}\enspace
Without loss of generality, we may assume that the sequence $\{p_n\}_{n=1}^{\infty}$ is nondecreasing. The hypothesis on $\psi$ implies that, for positive integers $n$,
\begin{equation}\label{theta_2n}
\psi(2n) - \psi(n) = O(1).
\end{equation}
Since the $f_n$ are trigonometric polynomials,
\begin{equation*}
\max\limits_{x \in \mathbb{T}} S^*(f_n, x) =: B_n < \infty.
\end{equation*}

Construct a subsequence of positive integers $\{n_k\}_{k=1}^{\infty}$. Let $n_1 = 1$, and for every integer $k>1$ choose $n_k$ in terms of $n_1 < n_2 < \ldots < n_{k-1}$ so that the following requirements are satisfied:
\begin{gather}
\sum\limits_{j=1}^{k-1} \frac{B_{n_j}}{\sqrt{A_{n_j}}} \leqslant \frac{\chi(0) \psi(0)}{2} \sqrt{A_{n_k}},\label{n_k_1}\\
p_{n_{k-1}} \leqslant \frac{\sqrt{A_{n_k}}}{2^k} \label{n_k_2}.
\end{gather}
Since condition~\ref{item0} of Lemma~\ref{lemma} implies~\eqref{m_sum} for the sequence of sets $\{E_{n_k}\}_{k=1}^{\infty}$, Theorem~\ref{Borel-Kantelli} yields a sequence $\{x_k\}_{k=1}^{\infty}$ such that
\begin{equation}\label{m_limsup}
mE := m\varlimsup\limits_{k \to \infty} E_{n_k}^{x_k} = 2 \pi.  
\end{equation}

Set
\begin{equation}\label{F}
F(x) := \sum\limits_{k=1}^{\infty} \frac{f_{n_k}(x-x_k)}{\psi(p_{n_k}) \sqrt{A_{n_k}}}.
\end{equation}
For any positive integer $k$, condition~\ref{item1} of Lemma~\ref{lemma} implies that the degrees of the polynomials $f_{n_j}$ are at most $2p_{n_k}$ for $j \leqslant k$. Thus, using condition~\ref{item2} of Lemma~\ref{lemma} and inequality~\eqref{n_k_2}, we obtain, for $m \geqslant 2p_{n_k}$,
\begin{gather}\label{E_n(F)}
E_m(F)_1 \leqslant \left\| \sum\limits_{j=k+1}^{\infty} \frac{f_{n_j}}{\psi(p_{n_j}) \sqrt{A_{n_j}}} \right\|_1  \leqslant C \sum\limits_{j=k+1}^{\infty} \frac{1}{\psi(p_{n_j}) 2^j}.
\end{gather}
Since $E_m(F)_1$ is nonincreasing,
\begin{equation}\label{E_n(F)_2}
\sum\limits_{m=2p_1}^{\infty} [\psi(m+1) - \psi(m)] E_m(F)_1 \leqslant \sum\limits_{k=1}^{\infty} [\psi(2p_{n_{k+1}}) - \psi(2p_{n_k})] E_{2p_{n_k}}(F)_1.
\end{equation}
Inequalities~\eqref{E_n(F)_2},~\eqref{E_n(F)}, and relation~\eqref{theta_2n} imply that there exists a constant $C'>0$ such that
\begin{multline*}
\sum\limits_{m=2p_1}^{\infty} [\psi(m+1) - \psi(m)] E_m(F)_1 \leqslant C \sum\limits_{k=1}^{\infty} \sum\limits_{j=k+1}^{\infty} \frac{\psi(2p_{n_{k+1}}) - \psi(2p_{n_k})}{\psi(p_{n_j}) 2^j} \\
= C \sum\limits_{j=2}^{\infty} \sum\limits_{k=1}^{j-1}  \frac{\psi(2p_{n_{k+1}}) - \psi(2p_{n_k})}{\psi(p_{n_j}) 2^j} \leqslant C \sum\limits_{j=2}^{\infty} \frac{\psi(2p_{n_j})}{\psi(p_{n_j}) 2^j} \leqslant C' \sum\limits_{j=2}^{\infty}  \frac{1}{2^j} < \infty.
\end{multline*}
Consequently, by Lemma~\ref{lemma_series}, the function $F$ satisfies condition~\eqref{series_cond}.

By~\eqref{n_k_2}, the series~\eqref{F} can be integrated term by term, and the Fourier coefficients of $F$ are equal to the sums of the Fourier coefficients of the functions
\begin{equation*}
\frac{f_{n_k}(x-x_k)}{\psi(p_{n_k}) \sqrt{A_{n_k}}},\quad k \in \mathbb{N}.
\end{equation*}
Therefore, the partial sums of the Fourier series of $F$ can be written as follows:
\begin{equation*}
S_m(F, x) = \sum\limits_{k=1}^{\infty} \frac{S_m(f_{n_k},x-x_k)}{\psi(p_{n_k}) \sqrt{A_{n_k}}},\quad m \in \mathbb{Z}_+.
\end{equation*}

Fix an integer $k>1$. Let $x \in E_{n_k}^{x_k}$. Since $x - x_k \in E_{n_k}$, condition~\eqref{S>m} of Lemma~\ref{lemma} implies that there exists $m_k = m_k(x) \leqslant p_{n_k}$ such that
\begin{equation}\label{S>m3}
|S_{m_k}(f_{n_k}, x-x_k)| > \psi(p_{n_k}) \chi(m_k) A_{n_k}.
\end{equation}
Consider the $m_k$th partial sum of the Fourier series of $F$ at $x$:
\begin{multline}\label{S_n}
S_{m_k}(F, x) = \sum\limits_{j=1}^{k-1} \frac{S_{m_k}(f_{n_j},x-x_j)}{\psi(p_{n_j}) \sqrt{A_{n_j}}} + \frac{S_{m_k}(f_{n_k},x-x_k)}{\psi(p_{n_k}) \sqrt{A_{n_k}}} \\
+ \sum\limits_{j=k+1}^{\infty} \frac{S_{m_k}(f_{n_j},x-x_j)}{\psi(p_{n_j}) \sqrt{A_{n_j}}} =: I_1 + I_2 + I_3.
\end{multline}
Using condition~\ref{item1} of Lemma~\ref{lemma} and property~\eqref{n_k_1}, we estimate
\begin{equation}\label{I1}
|I_1| \leqslant \sum\limits_{j=1}^{k-1} \frac{B_{n_j}}{\psi(0) \sqrt{A_{n_j}}}  \leqslant \frac{\chi(0)}{2} \sqrt{A_{n_k}}.
\end{equation}
Applying the estimate for the Dirichlet kernel, condition~\ref{item2} of Lemma~\ref{lemma}, and condition~\eqref{n_k_2}, we obtain
\begin{equation}\label{I3}
|I_3| \leqslant C \sum\limits_{j=k+1}^{\infty} \frac{p_{n_k} + 1/2}{\psi(0)\sqrt{A_{n_j}}} \leqslant \frac{2C}{\psi(0)} \sum\limits_{j=k+1}^{\infty} \frac{p_{n_k}}{p_{n_{j-1}} 2^j} \leqslant \frac{2C}{\psi(0)}.
\end{equation}
Finally, inequality~\eqref{S>m3} yields
\begin{equation*}
|I_2| > \chi(m_k) \sqrt{A_{n_k}}.
\end{equation*}
From this, representation~\eqref{S_n}, and inequalities~\eqref{I1} and~\eqref{I3}, we obtain, for $x \in E_{n_k}^{x_k}$,
\begin{equation}\label{S_n>}
\max\limits_{0 \leqslant m \leqslant p_{n_k}} \frac{|S_m(F, x)|}{\chi(m)} \geqslant \frac{|I_2|}{\chi(m_k)} - \frac{|I_1|}{\chi(0)} - \frac{|I_3|}{\chi(0)} > \frac{1}{2} \sqrt{A_{n_k}} - \frac{2C}{\chi(0) \psi(0)}.
\end{equation}
By~\eqref{m_limsup}, the set $E$ has full measure. For every $x \in E$, we have $x \in E_{n_k}^{x_k}$ for infinitely many $k \in \mathbb{N}$. Inequality~\eqref{S_n>} then implies~\eqref{not_est}. This proves Lemma~\ref{lemma2}.

\section{Proofs of the main results}

\subsection{Proof of Theorem~\ref{theorem1}}
First suppose that condition~\ref{item:first} of Theorem~\ref{theorem1} holds, and let us show that condition~\ref{item:second} holds. Fix a function $f \in L(\mathbb{T})$ satisfying~\eqref{series_cond} that is not a trigonometric polynomial. Let $n_0 = 0$, and for $k \in \mathbb{N}$ define
\begin{equation}\label{n_k}
n_k := \min\left\{n \in \mathbb{Z}_+: E_n(f)_1 < 2^{-k} E_0(f)_1 \right\}.
\end{equation}
It follows from definition~\eqref{n_k} that
\begin{equation}\label{E_0(f)}
2^{-(k+1)} E_0(f)_1 \leqslant E_n(f)_1 < 2^{-k} E_0(f)_1,\quad n \in \{n_k, \dots,  n_{k+1}-1\}.
\end{equation}

By estimates~\eqref{E_0(f)}, we have
\begin{multline}\label{sum_n_k}
\sum\limits_{n=1}^{\infty} [\psi(n+1) - \psi(n)] E_n(f)_1 \geqslant \sum\limits_{k=0}^{\infty} \sum\limits_{n = n_k}^{n_{k+1}-1} \frac{1}{2^{k+1}}[\psi(n+1) - \psi(n)] E_0(f)_1 \\
= \frac{1}{4} \sum\limits_{k=0}^{\infty}[\psi(n_{k+1}) - \psi(n_k)] E_0(f)_1 \sum\limits_{l=k}^{\infty} \frac{1}{2^l} = \frac{1}{4} \sum\limits_{l=0}^{\infty} \frac{1}{2^l} \sum\limits_{k=0}^{l} [\psi(n_{k+1}) - \psi(n_k)] E_0(f)_1 \\
= \frac{1}{4} \sum\limits_{l=0}^{\infty} \frac{[\psi(n_{l+1}) - \psi(n_0)] E_0(f)_1}{2^l} > \frac{1}{4} \sum\limits_{l=0}^{\infty} E_{n_l}(f)_1 \psi(n_{l+1}) - \frac{\psi(n_1) E_0(f)_1}{2}.
\end{multline}
Since, by Remark~\ref{remark3}, the series $\sum\limits_{n=1}^{\infty} [\psi(n+1) - \psi(n)] E_n(f)_1$ converges, estimate~\eqref{sum_n_k} implies that the series $\sum\limits_{k=0}^{\infty} E_{n_k}(f)_1 \psi(n_{k+1})$ converges as well. Hence, there exists a sequence of positive numbers $\{\alpha_k\}_{k=1}^{\infty}$ such that $\alpha_k \to \infty$ as $k \to \infty$ and
\begin{equation}\label{alpha_series}
\sum\limits_{k=0}^{\infty} \alpha_k E_{n_k}(f)_1 \psi(n_{k+1}) < \infty.
\end{equation}
Define $\theta(m) := \alpha_k$ for $m \in \{n_k, \ldots, n_{k+1}-1\}$.

For $y > 0$ and $k \in \mathbb{Z}_+$, consider the set
\begin{equation}\label{G_n}
G_k(y) := \left\{x \in \mathbb{T}: \max\limits_{n_k \leqslant m < n_{k+1}} \frac{\theta(m) |f(x) - S_m(f, x)|}{\chi(m)} > y \right\}.
\end{equation}
Fix $k \in \mathbb{Z}_+$ such that $n_{k+1} > n_k$, and estimate $m G_k(y)$. Let $T_k$ be a polynomial of best approximation of degree at most $n_k$ to $f$ in $L(\mathbb{T})$, and put $g_k(x) := f(x) - T_k(x)$. Then
\begin{equation}\label{|g_k|}
\| g_k \|_1 = E_{n_k}(f)_1,\quad k \in \mathbb{N}.
\end{equation}
Chebyshev's inequality implies that
\begin{equation}\label{Cheb}
m\left\{x \in \mathbb{T}: |g_k(x)| > y \right\} \leqslant \frac{\|g_k\|_1}{y}.
\end{equation}
Since the sequence of sublinear operators $\{M_n(\psi, \chi)\}_{n \in \mathbb{Z}_+}$ satisfies a weak-type $(1, 1)$ inequality, inequality~\eqref{Cheb} implies that there exists a constant $C_1>0$ such that, for every $f \in L(\mathbb{T})$, $k \in \mathbb{N}$, and $y>0$,
\begin{equation}\label{weak_n2}
m\left\{x \in \mathbb{T}: \max\limits_{n_k \leqslant m < n_{k+1}} \frac{|S_m(g_k, x)|}{\chi(m)} + |g_k(x)| > y \right\} \leqslant \frac{C_1 \psi(n_{k+1})}{y} \|g_k\|_1.
\end{equation}
Since the degree of $T_k$ is at most $n_k$, we obtain
\begin{equation*}
\max\limits_{n_k \leqslant m < n_{k+1}} |f(x) - S_m(f, x)| \leqslant \max\limits_{n_k \leqslant m < n_{k+1}} |S_m(g_k, x)| + |g_k(x)|.
\end{equation*} 
This, inequality~\eqref{weak_n2}, and identity~\eqref{|g_k|} yield
\begin{equation}\label{mG_k}
m G_k(y) \leqslant \frac{C_1}{y} \alpha_k E_{n_k}(f)_1 \psi(n_{k+1}),\quad k \in \mathbb{N}.
\end{equation}
Then estimate~\eqref{mG_k} and property~\eqref{alpha_series} imply that there exists a constant $C_2$ such that, for every $y>0$,
\begin{equation*}
m\left\{x \in \mathbb{T}: \sup\limits_{m \in \mathbb{Z}_+} \frac{\theta(m)}{\chi(m)} |f(x) - S_m(f, x)| > y \right\} \leqslant \sum\limits_{k=0}^{\infty} m G_k(y) \leqslant \frac{C_2}{y}.
\end{equation*}
Consequently, the quantity
\begin{equation*}
\varlimsup\limits_{m \to \infty} \frac{\theta(m) |f(x) - S_m(f, x)|}{\chi(m)}
\end{equation*}
is finite almost everywhere, and thus condition~\ref{item:second} of Theorem~\ref{theorem1} holds.

Now suppose that condition~\ref{item:first} of Theorem~\ref{theorem1} fails. Lemmas~\ref{lemma} and~\ref{lemma2} imply the existence of a function $F \in L(\mathbb{T})$ satisfying~\eqref{series_cond} for which relation~\eqref{not_est} holds almost everywhere. Consequently, condition~\ref{item:second} of Theorem~\ref{theorem1} also fails. This proves Theorem~\ref{theorem1}.

\begin{remark}
In the proof of the implication from condition~\ref{item:first} to condition~\ref{item:second}
in Theorem~\ref{theorem1}, the monotonicity of $\chi$ was not used.
Thus, for this implication, $\chi$ may be taken to be an arbitrary positive
function. This makes it possible, when $\chi$ tends to zero, to establish a rate of convergence of the partial sums to a function in $H_1^{\omega}$ on a set of full measure. Moreover, if condition~\eqref{series_cond} is replaced from the outset by
condition~\eqref{series_E_n(f)}, the appeal to Remark~\ref{remark3} becomes unnecessary, and the requirement
that the sequence $n(\psi(n+1)-\psi(n))$ be nonincreasing can be dropped. Thus, the following more general statement holds: if $\psi$ is positive and nondecreasing, $\chi$ is an arbitrary positive function, the sequence $\{M_n(\psi,\chi)\}_{n\in\mathbb Z_+}$ satisfies a weak-type $(1,1)$ inequality, and $f\in L(\mathbb T)$ satisfies condition~\eqref{series_E_n(f)}, then relation~\eqref{max_S_n} holds for almost every $x\in\mathbb T$.
\end{remark}

\subsection{Proof of Corollary~\ref{theorem2}}

To prove Corollary~\ref{theorem2}, we need a lemma from~\cite[Lemma~3.5]{Konyagin}. We state it in the form needed here.

\begin{thmx}
Let $n$ be a sufficiently large positive integer, and suppose that
\begin{equation}\label{psi_lnln}
\theta(n) = o\left(\sqrt{\frac{\ln{n}}{\ln{\ln{n}}}}\right)
\end{equation}
as $n \to \infty$. Then there exists a nonnegative trigonometric polynomial $T_n$ of degree at most $n$ such that $\|T_n\|_1 = 2\pi$ and, for every $x \in \mathbb{T}$,
\begin{equation}
\max\limits_{0 \leqslant m \leqslant n} S_m(T_n, x) > \theta(n) A_n,
\end{equation} 
where $A_n \to \infty$ as $n \to \infty$.
\end{thmx}

Put $\theta=\psi\chi$. By hypothesis, $\theta(n)=o\bigl(\sqrt{\ln n/\ln\ln n}\bigr)$, so Konyagin's lemma gives a nonnegative trigonometric polynomial $T_n$ of degree at most $n$ with $\|T_n\|_1=2\pi$ and
\begin{equation*}
\max_{0\le m\le n}S_m(T_n,x)>\psi(n)\chi(n)A_n,
\end{equation*}
where $A_n\to\infty$. Set $p_n=n$, $f_n=T_n$, and $E_n=\mathbb{T}$. Then the conditions of Lemma~\ref{lemma} are satisfied, and the required inequality~\eqref{S>m} follows from the monotonicity of $\chi$ and the nonnegativity of $S_m(T_n,x)$. Applying Lemma~\ref{lemma2}, we obtain the required function $F$.

\end{document}